\input amstex.tex
\input amsppt.sty   
\magnification 1200
\vsize = 9.5 true in
\hsize=6.2 true in
\NoRunningHeads        
\parskip=\medskipamount
        \lineskip=2pt\baselineskip=18pt\lineskiplimit=0pt
       
        \TagsOnRight
        \NoBlackBoxes

        \topmatter
        \title
        Bounded Sobolev Norms for Linear Schr\"odinger Equations 
        Under Resonant Perturbations       
        \endtitle
        \author
        W.-M.~Wang
        \endauthor
\address
Departement de Mathematique, Universite Paris Sud, 91405 Orsay Cedex, FRANCE
\endaddress
        \email
{wei-min.wang\@math.u-psud.fr}
\endemail
\thanks  
{I thank Fudan University for hospitality, where part of this work was done.}
\endthanks

        \bigskip\bigskip
        \bigskip
        \toc
        \bigskip
        \bigskip 
        \widestnumber\head {Table of Contents}
        \head 1. Introduction and statement of the theorems\endhead
        \head 2. Pure point spectrum for the Floquet Hamiltonian 
        \endhead
        \head 3. A Newton scheme
        \endhead
        \head 4. Computations of local eigenvalues and eigenfunctions
        \endhead
        \head 5. Proof of the theorems
        \endhead
        \endtoc
        \endtopmatter
        \vfill\eject
        \bigskip
\document
\head{\bf 1. Introduction and statement of the theorems}\endhead
We consider time periodic Schr\"odinger equations with periodic boundary conditions
of the form 
$$i\frac\partial{\partial t}u=-\Delta u+V(x, t)u,
\tag 1.1$$
where $x\in\Bbb T$, $t\in\Bbb R$, and $V$ is a real analytic potential periodic in $x$ and $t$.
The spectrum of the Laplacian: $\sigma(-\Delta)=\{j^2, j\in\Bbb Z\}$. We are interested in 
{\it resonant} perturbations, i.e., when the frequency $\omega$ of the time periodic potential
$V$ is an integer, $\omega\in\Bbb Z\backslash\{0\}$.

Assume $u$ is a solution, the $L^2$ norm is conserved by the Schr\"odinger flow map:
$$\Vert u(t)\Vert_{L^2(\Bbb T)}=\Vert u(0)\Vert_{L^2(\Bbb T)},\tag 1.2$$
for all $t$, and if $u(0)\in H^s(\Bbb T)$, ($s>0$), then $u(t)\in H^s(\Bbb T)$ for all $t$.
Unlike the time independent Schr\"odinger equations, in general, the $H^s$ norms 
of solutions to time dependent equations as in (1.1) can grow in time. 
Here we are concerned with the bounds on the Sobolev norms:
$\Vert u(t)\Vert_{H^s(\Bbb T)}$ as $t\to\infty$ when $V$ is resonant. This is generally 
speaking a more ``dangerous" case, where there is possible growth of $H^s$ norms.

We note that for a linear equation of the form (1.1), if one assumes $V$ is smooth in $x$
and $t$ (not necessarily periodic), one has the a priori bound (cf. \cite{B2, 3}):
$$\Vert u(t)\Vert_{H^s}\leq C_s(1+|t|^s)\Vert u(0)\Vert_{H^s},\tag 1.3$$
where $u(t)$ is the solution to (1.1) with initial condition $u(0)$.

Under a natural spectral condition (cf. (H1) of Theorem 0 in sect. 5), we prove that the $H^s$
norms of solutions to (1.1) remain bounded for all $s>0$, 
$$\Vert u(t)\Vert_{H^s}\leq C_s\Vert u(0)\Vert_{H^s},\tag 1.4$$
for all $t$. We show that this spectral condition (H1) is verified for small potentials $V$. 

Previously in \cite{B1}, it was shown that for time quasi-periodic potentials with Diophantine
frequencies (hence non-resonant)
$$\Vert u(t)\Vert_{H^s}\leq C_s(\log (1+|t|))^{C_s}\Vert u(0)\Vert_{H^s},\tag 1.5$$
for the corresponding solutions to (1.1). (1.5) holds in 1-d and 2-d when the time quasi-periodic
potential $V$ is small. In the periodic case, (1.5) was observed by T. Spencer \cite{S}, with 
no assumptions on the frequency $\omega$. The present paper constructs an explicit
example where there is {\it no} growth of Sobolev norms. It is partially motivated by results in
\cite {B1-3, S}.

In a companion paper \cite{W}, using related constructions, we show that for a general bounded,
time dependent potential $V(x,t)$, $x\in\Bbb T$ and $t\in\Bbb R$, which is analytic and periodic in $x$, smooth in $t$ (with no further specifications on the time dependence), the growth of Sobolev 
norms is at most logarithmic in $t$. Previously, Bourgain \cite{B2,3} showed that the 
growth of Sobolev norms is at most polynomial in $t$: $t^\epsilon$ (for any $\epsilon>0$) for $V(x,t)$ bounded, periodic in $x$, smooth in $x$ and $t$.
\smallskip
\noindent {\it The Floquet Hamiltonian}

When $V$ is periodic in time, it is well known from \cite{EV, H, YK} that properties of the solutions
to (1.1) can be deduced from the spectral properties of the corresponding Floquet Hamiltonian:
$$H=-i\frac{\partial}{\partial t}-\Delta+V\tag 1.6$$
on $L^2(\Bbb T_x\times\Bbb T_{t, \omega})$, where $\Bbb T_x=[0, 2\pi)$ with periodic boundary
conditions and $T_{t, \omega}=[0, 2\pi/\omega)$ with periodic boundary conditions.

By Fourier series, $H$ is unitarily equivalent to 
$$\hat H=\text { diag }(n\omega+j^2)+\hat V*\text   {  on }\ell^2(\Bbb Z^2),\tag 1.7 $$
where $\hat V(j, n)$ are the Fourier coefficients of $V$: 
$$V(x, t)=\sum_{(j,n)\in\Bbb Z^2}\hat V(j, n)e^{i(jx+n\omega t)},\tag 1.8$$
and $\hat V*$ denotes convolution:
$$[ \hat V*u](j, n)=\sum_{(j',n')}\hat V(j-j', n-n')u(j',n').\tag 1.9$$

In this paper, for simplicity, we take $$V(x, t)=2\cos x\cos t, \tag 1.10$$
instead of a more general analytic periodic potential. The frequency $\omega=1$ here. The method
here applies in the general case. (The term $2\hat V(0,1)\cos t$ can be eliminated by replacing $u$ by $ue^{2i\sin t\hat V(0,1)}$. This elimination procedure clearly holds more generally for potentials which only depend on $t$.) The Floquet Hamiltonian is then 
$$\hat H=\text{diag }(n+j^2)+\tilde\Delta\tag 1.11$$
where 
$$\aligned \tilde\Delta (j,n;j',n')&=1,\quad |j-j'|_{\ell^1}=1\text{ and }|n-n'|_{\ell^1}=1 ,\\
&=0,\quad \text{otherwise}.\endaligned\tag 1.12$$

Writing $H$ for $\hat H$ from now on and adding a parameter $\delta$ in front of $\tilde\Delta$
(as part of the arguments need $\delta$ to be small), in the rest of the paper, we shall study the 
spectral properties of the operator
$$H=\text{diag }(n+j^2)+\delta\tilde\Delta\text {  on }\ell^2(\Bbb Z^2)\tag 1.13$$
with $\tilde\Delta$ defined as in (1.12). For simplicity, we also denote 
$|\,|_{\ell^1}$ by $|\,|$.

It is known from \cite{EV} that $H$ has pure point spectrum for all $\delta$ by using compactness
arguments. However to bound the Sobolev norms, we need localization properties of the 
eigenfunctions of $H$. We have 

\proclaim{Theorem 1} There exist $0<K<\infty$, $0<c<\infty$ and $0<C<\infty$ such that
for $\delta$ small enough, the eigenfunctions $\phi$ of the Floquet Hamiltonian (1.11) with
eigenvalue $E$ satisfy either
$$\align &|\phi(x)|\leq Ce^{-\frac{|x-(0,[E])|}{K}}\tag 1.14\\
\text{or }&|\phi(x)|\leq C\sum_{i=1}^2e^{-c|x-x_i|}\tag 1.15\endalign$$
for some $x_i=(n_i,\pm j_i)$ satisfying $|n_i-E|\geq K$ and $|n_i+j_i^2-E|\leq 2\delta$, 
where $[E]$ is the integer part of $E$. 
\endproclaim

As a direct consequence of Theorem $1$, we have 

\proclaim{Theorem 2} Let $s>0$ and $u(0)\in H^s(\Bbb T)$. Then the solution $u(t)$ to 
(1.1, 1.10) with the initial condition $u(0)$ is in $H^s(\Bbb T)$ for all $t$ and satisfies
$$\Vert u(t)\Vert_{H^s}\leq C_s\Vert u(0)\Vert_{H^s}.\tag 1.16$$
\endproclaim

Theorems $1$ and $2$ will be proved in sect. 5. Theorem $2$ follows from (1.14, 1.15)
and a standard dyadic expansion. The main work is the proof of Theorem $1$. It is an 
Anderson localization (A. L.) type of results in the Fourier space. The main novelty
is that it holds for a {\it fixed} potential. In the usual A. L. setting, the potential depends 
on a {\it parameter} and localization holds on a set of paramerters with large or (some times)
full measure, cf. e.g., \cite {FS, GB, GK}. The reason we do not need a parameter here is because of the separation
properties of the set $\{j^2,\,j\in\Bbb Z\}$. We use this to prove spacing of local eigenvalues
and  then {\it uniformly} localized eigenfunctions.
\smallskip
\head{\bf 2. Pure point spectrum for the Floquet Hamiltonian}\endhead
From (1.13), when $V(x,t)=2\cos x\cos t$, the Floquet Hamiltonian $H$ has the form
$$H=\text{diag }(n+j^2)+\delta\tilde\Delta\text {  on }\ell^2(\Bbb Z^2),\tag 2.1$$
where 
$$\aligned \tilde\Delta (j,n;j',n')&=1,\quad |j-j'|=1\text{ and }|n-n'|=1 ,\\
&=0,\quad \text{otherwise}.\endaligned\tag 2.2$$ 
As mentioned in sect. 1, $H$ has pure point
spectrum \cite{EV}. For completeness we give a proof below to this fact using our formalism.
Since $\delta$ will be taken to be small in Theorems $1$ and $2$, we will only address that
case. (The general scheme presented in this paper extend to arbitrary $\delta$, although
some of the conditions, cf. Theorem $0$ in sect. 5, are verified only for small $\delta$ 
for the moment.)

\proclaim{Lemma 2.1} $H$ has pure point spectrum for $|\delta|<1/4$.
\endproclaim
\demo{Proof} When $\delta=0$, $\sigma(H)=\{n+j^2\}=\Bbb Z$ with infinite multiplicity.
Let $P$ be the projection onto the eigenspace of eigenvalue $0$ and $P^c$ the projection
onto the complement. $P^c=\oplus_{N\in\Bbb Z\backslash\{0\}}P_N$, where $P_N$ is the 
projection onto the eigenspace of eigenvalue $N$, $N\neq 0$. When $\delta\neq 0$ ($0<\delta<1/4$),
$\sigma(H)\subseteq N+[-2\delta, 2\delta]$, $N\in\Bbb Z$. It is sufficient to look at spectral 
parameters $E$ such that $|E|\leq 2\delta$. 

Using the Feshbach projection or equivalently Grushin problem method, cf. e.g., \cite{SZ}, we 
have that $E\in\sigma(H)$ if and only if $0\in\sigma(\tilde H_E)$, where 
$$\tilde H_E=E-H^{00}-H^{0c}(E-H^{cc})^{-1}H^{c0}\tag 2.3$$
and 
$$\aligned &H^{00}=PHP=\delta\bar\Delta\\
&\bar\Delta(j,n;j'n')=\cases\tilde\Delta(j,n;j'n')\quad (j,n),\,(j',n')\in\{(0,0),(1,-1),(-1,-1)\}\\
0\text{ otherwise.}\endcases\\
&H^{0c}=PHP^c=\delta P\tilde\Delta P^c\\
&H^{c0}=P^cHP=\delta P^c\tilde\Delta P\\
&H^{cc}=P^cHP^c\endaligned\tag 2.4$$

We note that $\bar\Delta$ is rank $3$. Let $\tilde\Delta_{0N}=P\tilde\Delta P_N$, ($N\neq 0$), then from
(2.2), $\tilde\Delta_{0N}(j,n;j'n')\neq 0$ only if 
$$N=n\pm 1+(j\pm 1)^2|_{n+j^2=0}=\pm 2j+\cases 2\\0,\endcases\tag 2.5$$
and $|j-j'|=1$ and $|n-n'|=1$. We have 
$$H^{0c}=\delta\oplus_{N\neq 0}\tilde\Delta_{0N}.\tag 2.6$$
and $$H^{c0}=\delta\oplus_{N\neq 0}\tilde\Delta_{N0}.\tag 2.7$$

Define $$\aligned A_E&=H^{00}+H^{0c}(E-H^{cc})^{-1}H^{c0}\\
&{\overset\text{\rm def }\to=}H^{00}+\delta^2B_E\quad (|E|\leq 2\delta).\endaligned\tag 2.8$$
Let $D$ be the diagonal part of $H$ and $D_{NN}=P_NDP_N$, $N\neq 0$. Using (2.6, 2.7) and the resolvent equation, we have 
$$\aligned B_E=&\oplus_{N\neq 0}\tilde\Delta_{0N}(E-D_{NN})^{-1}\tilde\Delta_{N0}\\
&\qquad +[\oplus_{N\neq 0}\tilde\Delta_{0N}(E-D_{NN})^{-1}]\delta\tilde\Delta^{cc}
(E-H^{cc})^{-1}
[\oplus_{N'\neq 0}\tilde\Delta_{N'0}],\endaligned\tag 2.9$$
where $\tilde\Delta^{cc}=P^c\tilde\Delta P^c$. $B_E$ acts on the eigenspace $n+j^2=0$. 
$$\tilde\Delta_{0N}(E-D_{NN})^{-1}(j,n;j'n')\neq 0$$ only if $n+j^2=0$, $|j-j'|=1$ and $|n-n'|=1$. Using
(2.5), $$|\tilde\Delta_{0N}(E-D_{NN})^{-1}(j,n;j'n')|\leq \Cal O(1/|j|)\,(|j|\gg 1)$$ and since 
$$\Vert (E-H^{cc})^{-1}\Vert \leq 1-2\delta\quad (0<\delta<1/4),$$
$B_E$ is compact for all $E\in[-2\delta, 2\delta]$. Therefore $A_E$ is compact as $H^{00}$ 
is a rank $3$ operator, with $0$ the only possible accumulation point. So for all $E\in[-2\delta, 2\delta]$,
$\tilde H_E$ has pure point spectrum with $E$ the only possible accumulation point. 
Clearly the above argument goes through for all $E\in N+[-2\delta, 2\delta]$, $N\in\Bbb Z$, with
$E-N$ replacing $E$. 

Coming back to $H$, $E\in\sigma(H)$ if and only if $0\in\sigma(\tilde H_E)$. This implies
that $H$ has pure point spectrum with $\Bbb Z$ the only possible accumulation points. Here we 
also used the fact that the reduction in (2.3) preserves spectral multiplicity (cf. \cite{SZ}).\hfill $\square$
\enddemo

As explained in sect. 1, in order to prove boundedness of Sobolev norms, we need to have precise 
localization properties of the eigenfunctions of $H$. For that purpose, it is essential 
to exhibit eigenvalue spacing. As earlier, we only need to look at $\sigma(H)\cap [-2\delta, 2\delta]$
($0<\delta<1/4$), as the eigenfunctions for the other intervals are just translates in the $n$
direction.

When $\delta=0$, $\sigma(H)=\{n+j^2\}=\Bbb Z$. From perturbation theory, the only 
equi-energy parabola of relevance for the spectral range $[-2\delta, 2\delta]$ is $n+j^2=0$ 
($0<\delta<1/4$). Using a Newton scheme, we compute the perturbed local eigenvalues.
The result gives the necessary eigenvalue spacing in order to prove localization of eigenfunctions.
\smallskip 

\head{\bf 3. A Newton scheme}\endhead

Let $H$ be a linear operator on $\ell^2(\Lambda)$, $\Lambda\subseteq\Bbb Z^d$. We write 
$$H=D+\Delta H,$$
where $D$ is diagonal. Without loss of generality, we may assume $(\Delta H)_{ii}=0$, for all 
$i\in\Lambda$. Let $E$ be an eigenvalue of $D$, then $E=D_{ii}$ for some $i\in\Lambda$,
with eigenfunction $u=\delta_i$. We assume 
$$\aligned &D_{ii}\neq D_{jj}\quad \forall j\neq i\\
\text {and }&\Vert \Delta H\Vert<\frac{1}{2}\inf_{j\neq i}|D_{ii}-D_{jj}|.\endaligned\tag 3.1$$
We call $i$ the {\it resonant} site. Let $\Cal R=\{i\}$. $\Cal R^c=\Lambda\backslash \{i\}$.
We compute the eigenvalues and eigenfunctions using the following iteration scheme.

\noindent{\it Remark.}
Under the assumption (3.1), $i$ is the only resonant site, cf. (3.9). The scheme below however
can be applied to cases where there is symmetry, (3.1) is violated and there are more than $1$ 
resonant site. For example it can be used to compute eigenvalue splitting for the (time independent)periodic Schr\"odinger operator in 1-d. We leave this aspect of things to a future publication. The approach here is different from the Raleigh-Schr\"odinger scheme in quantum mechanics. 
It is closer to the Grushin-Feshbach effective operator method. More precisely it provides a way to compute eigenvalues when the effective operator is finite dimensional.

We seek solutions to the eigenvalue problem
$$(H-E)u=0,\tag 3.2$$
such that $u|_{\Cal R}=1$ is fixed.
As a zeroth order approximation 
$$E=D_{ii},\,u=\delta_i.\tag 3.3$$
So $$(H-E)u=\Delta H\delta_i{\overset\text{\rm def }\to=}F(u),\tag 3.4$$
where $F(u)$ is the error satisfying
$$F(u)|_{\Cal R}=0, F(u)=F(u)|_{\Cal R^c}.$$ 

Assume we have (3.2, 3.3) at the 
$n^{\text {th}}$ iteration, with $u=u^{(n)}$, $E=E^{(n)}$. To obtain the $(n+1)$th approximant, we 
write 
$$\aligned &u^{(n+1)}=u^{(n)}+\Delta u^{(n+1)}{\overset\text{\rm def }\to=}u+\Delta u,\\
&E^{(n+1)}=E^{(n)}+\Delta E^{(n+1)}{\overset\text{\rm def }\to=}E+\Delta E,\endaligned\tag 3.5$$
such that 
$$\align &(H-E^{(n)})u^{(n+1)}|_{{\Cal R}^c}=(H-E)(u+\Delta u)|_{{\Cal R}^c}=0,\tag 3.6\\
&(H-E^{(n+1)})u^{(n+1)}|_{\Cal R}=(H-E-\Delta E)(u+\Delta u)|_{\Cal R}=0\tag 3.7\endalign$$
are verified.

Since $\Delta u|_{\Cal R}=0$, from (3.6), $$(H-E)|_{{\Cal R}^c}\Delta u|_{{\Cal R}^c}=-(H-E)u|_{{\Cal R}^c}=-F(u)|_{{\Cal R}^c},\tag 3.8$$
so 
$$\Delta u|_{{\Cal R}^c}=-(H|_{{\Cal R}^c}-E)^{-1}F(u)|_{{\Cal R}^c}.\tag 3.9$$

From (3.7), 
$$(H-E)u|_{\Cal R}+(H-E)\Delta u|_{\Cal R}-\Delta E u|_{\Cal R}-\Delta E\Delta u|_{\Cal R}=0.$$
The first term is $0$ from (3.7), the fourth term is $0$ since $\Delta u|_{\Cal R}=0$. So
$$\Delta E=(H-E)\Delta u|_{\Cal R}=(\Delta H\Delta u)|_{\Cal R},\tag 3.10$$
where we used $u|_{\Cal R}=1$ by definition.
The error of approximation
$$\aligned F(u+\Delta u)&=(H-E-\Delta E)(u+\Delta u)\\
&=(H-E-\Delta E)(u+\Delta u)|_{{\Cal R^c}}\text { from } (3.7)\\
&=-\Delta E(u+\Delta u)|_{{\Cal R^c}}.\endaligned\tag 3.11$$
We now show that the above iteration scheme converges for $H$ in (2.1) restricted
to appropriate subsets of $\Bbb Z^2$.

\noindent{\it Convergence of the Newton scheme}

It is sufficient to look at $\sigma(H)\cap [-2\delta, 2\delta]$ ($0<\delta<1/4$). $n+j^2=0$ is the 
resonant parabola. Let $P=\{(j,n)|n+j^2=0\}$. $\Bbb Z^2\backslash P$ are non-resonant. For 
any two points $(j,n)$, $(j',n')\in P$, $(j,n)\neq(j',n')$, $|(j,n)-(j',n')|_\infty\geq d$ with   
$$\aligned d&=\max (||j'|-|j||(|j'|+|j|),\, |j-j'|)\\
&=\cases 2|j|,\quad |j|=|j'|, j\neq j'\\(|j'|-|j|)(|j'|+|j|),\quad |j|\neq |j'|.\endcases\endaligned\tag 3.12$$
For all $(j,n)\in P$, $|j|>1$, define 
$$\Lambda_j=\{(j', n')||(j',n')-(j,n)|_{\infty}\leq L_j\},\quad |j|\leq L_j\leq 2(|j|-1)\tag 3.13$$
to be the square centered at $(j,n)$ with side length $2L_j$. From (3.12), $P\cap\Lambda_j=\{(j,n)\},\,|j|>1$. So $(j,n)$ is the only resonant site in $\Lambda_j$ at $E=0$. 

For any $S\subset\Bbb Z^2$, define 
$$\aligned &H_S(j',n';j'',n'')=H(j',n';j'',n''),\quad (j',n'), (j'',n'')\in S\\
&H_S=0\text {  otherwise}.\endaligned \tag 3.14$$
We now prove that the Newton scheme in (3.2-3.11) converges for $H_{\Lambda_j}$, when 
$|j|$ is sufficiently large.

For simplicity of notation, we write $H$ for $H_{\Lambda_j}$. $\Cal R=\{(j,n)\}$, $(j,n)\in P$, $\Cal R^c=\Lambda_j\backslash\{(j,n)\}$. $H_{\Cal R^c}$ is $H$ restricted to $\Cal R^c$. Define 
$$F^{(k)}(u^{(k)}){\overset\text{\rm def }\to=}(H-E^{(k)})u^{(k)}.\tag 3.15$$
From (3.11)
$$F^{(k)}(u^{(k)})=-\Delta E^{(k)}u^{(k)}|_{\Cal R^c}.\tag 3.16$$
From (3.9, 3.10)
$$\align&\Delta u^{(k+1)}|_{\Cal R^c}=-(H_{\Cal R^c}-E^{(k)})^{-1}F^{(k)}(u^{(k)})|_{\Cal R^c},\tag 3.17\\
&\Delta E^{(k)}=\delta\sum_{\Sb |j-j'|=1\\|n-n'|=1\endSb}\Delta u^{(k)}(j',n'),\tag 3.18\endalign$$
where we have put back the superscript according to (3.5).

The following lemma shows that the Newton scheme (3.15-3.18) converges exponentially 
fast for $|j|$ sufficiently large.
\proclaim{Lemma 3.1} $$\align \Vert F^{(k)}\Vert &\leq \frac{2C^2\Vert F^{(0)}\Vert}{|j|}\Vert F^{(k-1)}\Vert<\frac{1}{2}\Vert F^{(k-1)}\Vert\text{ all } k\geq 1,\tag 3.19\\
|\Delta E^{(k)}|&<\frac{C}{|j|}\Vert F^{(k-1)}\Vert \text{ all } k\geq 1,\tag 3.20\endalign$$
where $F^{(k)}=F^{(k)}(u^{(k)})$, $1<C\leq 1/(1-2\delta)$, ($0<\delta<1/4$), $u^{(0)}=\delta_{(j,n)}$,
$E^{(0)}=0$, provided $|j|>4C^2\Vert F^{(0)}\Vert {\overset\text{\rm def }\to=}j_0$.
\endproclaim

\noindent{\it Remark.} (3.17-3.20) show that the above iteration scheme provides a convergent 
series expansion for the eigenvalue of $H_{\Lambda_j}$ in $[-2\delta, 2\delta]$ and its 
eigenfunction, although for the purpose of this paper, this is not needed. 

\demo{Proof} We start from $k=1$. From (3.17)
$$\Delta u^{(1)}|_{\Cal R^c}=-(H_{\Cal R^c}-E^{(0)})^{-1}F^{(0)}|_{\Cal R^c},\tag 3.21$$
so $$\Vert \Delta u^{(1)}\Vert \leq \Vert F^{(0)}\Vert,\tag 3.22$$
where we used $$\Vert H_{\Cal R^c}-E^{(0)}\Vert\geq 1.\tag 3.23$$
From (3.18),
$$\aligned \Delta E^{(1)}&=\delta\sum_{\Sb |j-j'|=1\\|n-n'|=1\endSb}\Delta u^{(1)}(j',n')\\
&=\delta\sum_{\Sb |j-j'|=1\\|n-n'|=1\endSb}[(D_{\Cal R^c}-E^{(0)})^{-1}F^{(0)}](j',n')\\
&\qquad\qquad+ [(D_{\Cal R^c}-E^{(0)})^{-1}\delta\tilde\Delta(H_{\Cal R^c}-E^{(0)})^{-1}F^{(0)}](j',n'),\endaligned\tag 3.24$$
so $$|\Delta E^{(1)}|\leq\frac{C}{|j|}\Vert F^{(0)}\Vert,\tag 3.25$$
where we used $$|n'+{j'}^2|=|n\pm 1+(j\pm 1)^2|\geq |j|,\, (|j|>1)\tag 3.26$$
and (3.23).

Using (3.25, 3.22) in (3.16), we have 
$$\Vert F^{(1)}\Vert \leq\frac{C}{|j|}\Vert F^{(0)}\Vert\cdot \Vert F^{(0)}\Vert.\tag 3.27$$
So $$\Vert F^{(1)}\Vert<\frac{1}{2}\Vert F^{(0)}\Vert,\tag 3.28$$
if $|j|>2C\Vert F^{(0)}\Vert$. (3.25, 3.28) show that (3.19, 3.20) hold at $k=1$. (3.25) shows
that $$\Vert H_{\Cal R^c}-E^{(1)}\Vert >1/C.\tag 3.29$$

Assume (3.19, 3.20) hold for all $k\leq K$, which implies also that 
$$\Vert H_{\Cal R^c}-E^{(K)}\Vert >1/C.\tag 3.30$$
From (3.17)
$$\Delta u^{(K+1)}|_{\Cal R^c}=-(H_{\Cal R^c}-E^{(K)})^{-1}F^{(K)}|_{\Cal R^c},\tag 3.31$$
so 
$$\Vert \Delta u^{(K+1)}|_{\Cal R^c}\Vert\leq C\Vert F^{(K)}\Vert\tag 3.32$$
using (3.30).  From (3.18)
$$\aligned \Delta E^{(K+1)}&=\delta\sum_{\Sb |j-j'|=1\\|n-n'|=1\endSb}\Delta u^{(K+1)}(j',n')\\
&=-\delta\sum_{\Sb |j-j'|=1\\|n-n'|=1\endSb}[(D_{\Cal R^c}-E^{(K)})^{-1}F^{(K)}](j',n')\\
&\qquad\qquad +[(D_{\Cal R^c}-E^{(K)})^{-1}\delta\tilde\Delta(H_{\Cal R^c}-E^{(K)})^{-1}F^{(K)}](j',n').\endaligned$$
So $$|\Delta E^{(K+1)}|\leq\frac{C}{|j|}\Vert F^{(K)}\Vert.\tag 3.33$$
Hence $$\Vert H_{\Cal R^c}-E^{(K+1)}\Vert >1/C.\tag 3.34$$

Using (3.33, 3.32, 3.17) in (3.16), we have 
$$\Vert F^{(K+1)}\Vert\leq\frac{C}{|j|}\Vert F^{(K)}\Vert\cdot C(\sum_{k=0}^K\Vert F^{(k)}\Vert),\tag 3.35$$
$$\Vert F^{(K+1)}\Vert\leq\frac{2C^2}{|j|}\Vert F^{(K)}\Vert\cdot \Vert F^{(0)}\Vert,\tag 3.36$$
if $|j|>4C^2\Vert F^{(0)}\Vert{\overset\text{\rm def }\to=}j_0$, where we used (3.19) from $1\leq k\leq K$
to estimate the sum in (3.35). (3.33, 3.34, 3.36) imply that the lemma holds by induction.\hfill $\square$
\enddemo
\smallskip
\head{\bf 4. Computation of local eigenvalues and eigenfunctions}\endhead
Let $P=\{(j,n)|n+j^2=0\}$ be the resonant parabola at $E=0$. Let $\Lambda_j$ be the square 
defined as in (3.13). We use the convergent Newton scheme to compute the eigenvalue 
$E\in\sigma(H_{\Lambda_j})\cap[-2\delta, 2\delta]$ ($0<\delta<1/4$).

Let $u^{(0)}=\delta_{(j,n)}$, $E^{(0)}=0$. Then 
$$F^{(0)}(j',n')=F^{(0)}(u^{(0)})(j',n')=\cases\delta,\quad |j-j'|=1, \text { and } |n-n'|=1\\0, \text{ otherwise.}\endcases\tag 4.1$$
Assume $|j|>4C^2\Vert F^{(0)}\Vert=j_0$, so that Lemma 3.1 is applicable.

\proclaim {Lemma 4.1} $E$ has the convergent series expansion 
$$E=-\frac{\delta^2}{j^2}+ \frac{a_4\delta^4}{j^4}+\cdots,\tag 4.2$$
where $a_4$ is independent of the $L_j$ in (3.13).
\endproclaim
\demo{Proof}
From (3.17),
$$\Delta u^{(1)}|_{\Cal R^c}=-(H_{\Cal R^c}-E^{(0)})^{-1}F^{(0)}|_{\Cal R^c}\tag 4.3$$
In view of (3.18), we compute $\Delta u^{(1)}$ on the set 
$$\{(j',n')||n-n'|=1\text{ and } |j-j'|=1\}.\tag 4.4$$ We have 
$$\aligned \Delta u^{(1)}(j',&n')=-(D-E^{(0)})^{-1}F^{(0)}(j',n')\\
&-(D-E^{(0)})^{-1}\delta\tilde\Delta(D-E^{(0)})^{-1}F^{(0)}(j',n')\\
&-(D-E^{(0)})^{-1}\delta\tilde\Delta(D-E^{(0)})^{-1}\delta\tilde\Delta(D-E^{(0)})^{-1}F^{(0)}(j',n')\\
&-(D-E^{(0)})^{-1}\delta\tilde\Delta(D-E^{(0)})^{-1}\delta\tilde\Delta(D-E^{(0)})^{-1}\delta\tilde\Delta(D-E^{(0)})^{-1}F^{(0)}(j',n')\\
&-(D-E^{(0)})^{-1}\delta\tilde\Delta(D-E^{(0)})^{-1}\delta\tilde\Delta(D-E^{(0)})^{-1}\delta\tilde\Delta\\
&\quad (D-E^{(0)})^{-1}\delta\tilde\Delta(D-E^{(0)})^{-1}F^{(0)}(j',n')\\
&-(D-E^{(0)})^{-1}\delta\tilde\Delta(D-E^{(0)})^{-1}\delta\tilde\Delta(D-E^{(0)})^{-1}\delta\tilde\Delta(D-E^{(0)})^{-1}\delta\tilde\Delta\\
&\quad(D-E^{(0)})^{-1}\delta\tilde\Delta(H_{\Cal R^c}-E^{(0)})^{-1}F^{(0)}(j',n').\endaligned\tag 4.5$$

In view of the support of $F^{(0)}$: $\text {supp }F^{(0)}$ and $\tilde\Delta$, the second and fourth 
terms in (4.5) are zero on the set in (4.4), the last two terms are of order $\Cal O(1/|j|^3)$. So 
$$\aligned \Delta u^{(1)}(j',n')=&- (D-E^{(0)})^{-1}F^{(0)}(j',n')\\
&-\delta^2(D-E^{(0)})^{-1}\tilde\Delta(D-E^{(0)})^{-1}\tilde\Delta(D-E^{(0)})^{-1}F^{(0)}(j',n')\\
&+\Cal O (1/|j|^3)\endaligned\tag 4.6$$
The second term is in fact of order $\Cal O (1/|j|^3)$. Assume $n'=n-1$, $j'=j-1$. The second term 
gives 
$$\aligned&(D-E^{(0)})^{-1}\tilde\Delta(D-E^{(0)})^{-1}\tilde\Delta(D-E^{(0)})^{-1}F^{(0)}(j',n')\\
=&\frac{1}{n-1+(j-1)^2}\frac{1}{n-2+j^2}(\frac{1}{n-1+(j-1)^2}+\frac{1}{n-1+(j+1)^2})+\Cal O(1/|j|^3)\\
=&\Cal O(1/|j|^3)\endaligned\tag 4.7$$
as claimed. Similar estimates holds for the other $(j',n')$ in (4.4). (3.18) gives
$$\aligned \Delta E^{(1)}=&\delta\sum_{\Sb |n'-n|=1\\|j'-j|=1\endSb}(D-E^{(0)})^{-1}F^{(0)}(j',n')+\Cal O(1/|j|^3)\\
&=\delta^2[\frac{1}{n+1+(j-1)^2}+\frac{1}{n+1+(j+1)^2}+\frac{1}{n-1+(j+1)^2}+\frac{1}{n-1+(j-1)^2}]\\
&\qquad +\Cal O(1/|j|^3)\\
&=-\frac{\delta^2}{j^2}+\Cal O(\frac{1}{|j|^3})\endaligned\tag 4.8$$
The $\Cal O(1/|j|^3)$ is in fact absent from symmetry arguments. Using (3.16, 3.17), (4.8) is the only 
contribution at order $\Cal O(1/j^2)$. So $E$ has the expansion in (4.2). \hfill $\square$
\enddemo
\smallskip
\head{\bf 5. Proof of the theorems}\endhead
Let $E$ be an eigenvalue of $H$ with eigenfunction $\phi$: 
$$(H-E)\phi=0,\tag 5.1$$
$\phi\in\ell^2(\Bbb Z^2)$ from Lemma 2.1. Write $x=(j,n)\in\Bbb Z^2$. We prove 
\proclaim{Theorem 0} Assume 
\item{(H1)} $\exists L_0\gg 1$, such that 
$\text{dist } (0,\sigma(H_{\Lambda_0})\backslash\{0\})\gg e^{-L_0}$, where $\Lambda_0=[-L_0^2+1, L_0^2-1]^2$,
\item{(H2)} $\exists\ell_0$, $0<\ell_0\leq L_0$ such that for $|j|>\ell_0$, Lemma 3.1 is available,
\item{(H3)} if $0\in\sigma(H_{\Lambda_0})$, then $0\in\sigma(H_{\Lambda})$ for all 
$\Lambda=[-L, L]^2\supset\Lambda_0$ with the same multiplicity.

\noindent Under the conditions (H1-3),  there exist $0<K<\infty$, $0<C<\infty$, $0<c<\infty$ (depending only
on $\delta$, $L_0$), such that either
$$\align&|\phi(x)|\leq Ce^{-\frac{|x-(0,E)|}{K}}\tag 5.2\\
\text{or }&|\phi(x)|\leq C\sum_{i=1}^2 e^{-c|x-x_i|}\tag 5.3\endalign$$
for some $x_i=(n_i,\pm j_i)$ satisfying 
$$|n_i-E|\geq K\text { and }|n_i+j_i^2-E|\leq 2\delta.\tag 5.4$$
\endproclaim

\noindent{\it Remark.} If (H1) holds for the cube $\Lambda_0=[-L_0^2+1, L_0^2-1]^2$, then it holds
for the cubes $\Lambda'=[-L_0^2+L, L_0^2-L]^2$, ($0\leq L\leq L_0$).
\demo{Proof of the theorem} As earlier, we may assume $E\in[-2\delta, 2\delta]$ ($0<\delta<1/4$)
without loss of generality. This is because if we define $\bar\phi$ to be 
$\bar\phi (\cdot, \cdot)=\phi(\cdot -N,\cdot),\,N\in\Bbb Z$. Then $(H-\bar E)\bar\phi=0$
with $\bar E=N+E$. So $\bar\phi$ has the same localization properties as $\phi$.

Let $$\Lambda=[-L^2+1, L^2-1]^2,\,L>L_0.\tag 5.5$$
So $\Lambda\supset\Lambda_0$. Let $P=\{(j,n)|n+j^2=0\}$ be the set of resonant points. 
Let $P_\Lambda=P\cap\Lambda$. We cover $P_\Lambda$ with $\Lambda_0$, $\Lambda_{\pm j}$
($L_0\leq j<L$), $\Lambda_{\pm j}$ as defined in (3.13). Denote by $\mu$ the eigenvalues 
of $H_{\Lambda_0}$ in $[-2\delta, 2\delta]$, and $\lambda_{\pm j}$ the eigenvalues of 
$H_{\Lambda_{\pm j}}$ in  $[-2\delta, 2\delta]$, i.e., $\{\mu\}=\sigma(H_{\Lambda_0})\cap[-2\delta, 2\delta]$, $\lambda_{\pm j}=\sigma(H_{\Lambda_{\pm j}})\cap[-2\delta, 2\delta]$, ($L_0\leq j<L$).

Below we show that uniformly in $\Lambda$, $\mu$, $\lambda_{\pm j}$ are approximate eigenvalues of
$H_\Lambda$ in $[-2\delta, 2\delta]$. (For precise meaning of this, see (5.16, 5.17).)  Therefore
as $\Lambda\nearrow\Bbb Z^2$, the nonzero but  ``small" eigenvalues of $H_\Lambda$
``come" from $H_{\Lambda_{\pm j}}$ ($L_0\leq j<L$) (see (5.18, 5.19)). 

Let $\Lambda_*$ denote either $\Lambda_0$ or $\Lambda_{\pm j}$ and
$$\partial\Lambda_*=\{(j',n')\in\Lambda_*|\exists (j'',n'')\in\Lambda\backslash\Lambda_*,\text{ such that }
|j'-j''|=1 \text { and } |n'-n''|=1\}\tag 5.6$$
be the interior boundary of $\Lambda_*$ relative to $\Lambda$. Let $P_{\Lambda_*}=P\cap\Lambda_*$.

Define $\Lambda'_*=\Lambda_*\backslash P_{\Lambda_*}$. We have 
$$\text{dist }(\sigma(H_{\Lambda'_*}), [-2\delta, 2\delta])\geq 1-2\delta\, (0<\delta<1/4).\tag 5.7$$
Let $H_{\Lambda'_*}$ be defined as in (3.14) and 
$$\Gamma=H_{\Lambda_*}-H_{\Lambda'_*}.\tag 5.8$$
Assume $\phi$ is the eigenfunction for the eigenvalue
$\lambda_*\in\sigma(H_{\Lambda_*})\cap [-2\delta, 2\delta]$:
$$(H_{\Lambda_*}-\lambda_*)\phi=0.\tag 5.9$$
We have the identity
$$\phi=-(H_{\Lambda'_*}-\lambda_*)^{-1}\Gamma\phi\tag 5.10$$
using (5.8). Using (5.6, 5.7), we then have for some $\alpha>0$ (depending only on $\delta$)
$$\aligned &\Vert\phi\Vert_{\ell^2(\partial\Lambda_0)}\leq e^{-\alpha L_0}\quad (L_0 \text { large enough })\\
\text {and } &\Vert\phi\Vert_{\ell^2(\partial\Lambda_j)}\leq e^{-\alpha |j|}.\endaligned\tag 5.11$$

Let $H_\Lambda$ be defined as in (3.14) and 
$$\tilde\Gamma=H_\Lambda-H_{\Lambda_*}.\tag 5.12$$
We compute 
$$\aligned &(H_\Lambda-\lambda_*)\phi\\
=&(H_{\Lambda_*}-\lambda_*)\phi+\tilde\Gamma\phi.\endaligned\tag 5.13$$
So 
$$\aligned \Vert (H_\Lambda-\lambda_*)\phi\Vert &\leq\Vert\tilde\Gamma\phi\Vert\\
&\leq\cases e^{-\alpha L_0}\text {   if } \Lambda_*=\Lambda_0,\\
e^{-\alpha |j|}\text {   if } \Lambda_*=\Lambda_{\pm j},\endcases\endaligned\tag 5.14$$
where we used (5.11). This shows that $\exists \lambda\in\sigma(H_\Lambda)\cap[-2\delta, 2\delta]$
such that 
$$\aligned &|\lambda_*-\lambda|\leq e^{-\alpha L_0}\text{   if } \lambda_*\in\sigma(H_{\Lambda_0})\cap[-2\delta, 2\delta]\\
\text{and }&|\lambda_*-\lambda|\leq e^{-\alpha |j|}\text{   if } \lambda_*\in\sigma(H_{\pm j})\cap[-2\delta, 2\delta].\endaligned\tag 5.15$$

Let $I=|P_\Lambda|$ be the number of resonant points in $\Lambda$. We further label the 
(normalized) local eigenfunctions in (5.9) as $\phi_1$, $\phi_2$,...$\phi_I$. Using (5.7, 5.10), 
the matrix $A$ with entries $A_{ij}=(\phi_i, \phi_j)$ is invertible. Combined with (5.15), this shows that 
$$\align &\lambda=\mu+\Cal O(e^{-\alpha L_0})\tag 5.16\\
\text{ or }& \lambda=\lambda_{\pm j}+\Cal O(e^{-\alpha |j|}),\tag 5.17\endalign$$
where $\mu$ are the eigenvalues of $H_{\Lambda_0}$ in $[-2\delta, 2\delta]$ and $\lambda_{\pm j}$
are the eigenvalues of $H_{\Lambda_{\pm j}}$ in $[-2\delta, 2\delta]$ and have the convergent
series expansion as in (4.2). 

Let $\mu_{\min}$ be the smallest (in absolute value) nonzero eigenvalue of $H_{\Lambda_0}$
in $[-2\delta, 2\delta]$,
$$K'=\delta/\sqrt{|\mu_{\min}|}>0\tag 5.18$$
and set $$\Lambda=[-{K'}^2+K', {K'}^2+K']^2.\tag 5.19$$
(H1,3,4.2,5.16,5.17) then imply
$$\text{dist } [\sigma(H_\Lambda)\cap[-2\delta, 2\delta], \sigma(H_{\Lambda_{\pm j}})\cap[-2\delta, 2\delta]]\geq \frac{\Cal O(\delta^2)}{|j|^3}\text{  for } j>K'.\tag 5.20$$

Given $E\in[-2\delta, 2\delta]$, we define
\item{$\bullet$} $E$ resonant with $\Lambda_{\pm  j}$, if 
$$|E-\lambda_{\pm  j}|\leq e^{-j^\rho}\,
(j>K', \, 0<\rho<1/2), \tag 5.21$$where $\Lambda_{\pm   j}$ are as defined in (3.13), $\lambda_{\pm  j}$
are the eigenvalues of $H_{\Lambda_{\pm  j}}$ in $[-2\delta, 2\delta]$ and have the convergent 
series expansion in (4.2).
\item{$\bullet$} $E$ resonant with $\Lambda$ (defined in (5.19)) if 
$$\text{dist }(E,\sigma(H_\Lambda))\leq e^{-{K'}^\rho}\tag 5.22$$
and $$|E-\lambda_{\pm j}|>e^{-{j}^\rho}\, (0<\rho<1/2),\, \forall j>K'.\tag 5.23$$

\noindent Otherwise we say that $E$ is non resonant. Since for all $E\in[-2\delta, 2\delta]$, either (5.21) holds for some $\Lambda_{\pm j}$ or (5.23) holds, if (5.23) is satisfied and (5.22) is not, then from (5.19), $E$ is not an eigenvalue
of $\sigma(H)$ as 
$$\text{dist }(\sigma(H), \sigma(H_\Lambda))\leq e^{-\Cal O(K')}\ll e^{-{K'}^\rho}\,(0<\rho<1/2).$$
So if $E$ is non resonant and $E\in\sigma(H)\cap[-2\delta, 2\delta]$, then $E$ is an accumulation
point since $\sigma(H)$ is pure point from Lemma 2.1. In this case $E=0$.

If $E$ is resonant with $\Lambda_j$ ($j>K'$), then 
$$|E|\geq\frac{\Cal O(\delta^2)}{j^2}\tag 5.24$$
and 
$$\text{dist }(E,\sigma(H_\Lambda)\backslash\{0\})\geq\frac{\Cal O(\delta^2)}{{K'}^3}\, (j>K')\tag 5.25$$
from (5.18, 5.19, 5.21), and for $|j'|\neq |j|$, $|j'|>K'$
$$\aligned \text{dist }(E, \sigma(H_{\Lambda_{\pm j'}}))&\geq |\frac{\Cal O(\delta^2)}{j^2}-\frac{\Cal O(\delta^2)}{{j'}^2}|-\Cal O(e^{-{j}^\rho})-\Cal O(e^{-{j'}^\rho})\\
&\geq\frac{\Cal O(\delta^2)}{[\min(|j|,|j'|)]^3},\endaligned\tag 5.26$$
where the $\Cal O(\delta^2)$ are uniform in $|j|,\,|j'|>K'$.
(5.24-5.26) together with (H1,3) imply that 
$$\text{dist }(E,\sigma(H_{{\Lambda'}\backslash\{(\pm j, n)\}}))\geq\frac{\Cal O(\delta^2)}{{|j|}^3}\,(|j|>K')\tag 5.27$$
for all $\Lambda'\subseteq \Bbb Z^2$, where the $\Cal O(\delta^2)$ is uniform in $\Lambda'$. 

If $E$ is resonant with $\Lambda$, then since 
$$|E-\lambda_{\pm j}|>e^{-{j}^\rho}\, (0<\rho<1/2),\, \forall j>K',\tag 5.28$$
by definition,  
$$\text{dist }(E,\sigma(H_{\Lambda'\backslash\Lambda}))\geq\Cal O(e^{-J^\rho}),\,\forall\,\Lambda'=-[J^2,J^2]^2,\,J>K'.\tag 5.29$$

Assume $E$ is resonant with $\Lambda_{\pm j}$. Let 
$$\Lambda'(k)=(j,n)+[-k,k]^2,\,j>0,\, n+j^2=0,\,k\geq 2(j-1)\tag 5.30$$
be cubes centered at $(j,n)$ with side lengths $2k$. Let $P=\{(j',n')|n'+{j'}^2=0\}$ as before.
Assume $\Lambda'(k)$ is such that
$$\text{dist }[(j',n'),\partial\Lambda'(k)]\geq |j'|\tag 5.31$$
for all $(j',n')\in P\cap\Lambda'(k)$, where $\partial\Lambda'(k)$ is the boundary of
$\Lambda'(k)$ defined as in (5.6). Assume also that if $\Lambda'(k)\cap\Lambda\neq\emptyset$,
where $\Lambda$ is defined in (5.18, 5.19), then $\Lambda'(k)\supset\Lambda$. Similarly
we define $$\Lambda''(k)=(-j,n)+[-k,k]^2,\, j^2+n=0\tag 5.32$$

Let $\Cal R=\{(j,n),(-j,n)\}$, $j^2+n=0$. Given any $x\in\Bbb Z^2$, it is easy to see that there exists $k$, such that either $x\in\Lambda'(k)$ or $x\in\Lambda''(k)$ or both, and 
$$\align 10\text{ dist }(x,\Cal R\cap\Lambda'(k))&\geq\text{dist }(x,\partial\Lambda'(k))\geq\text{dist }
(x,\Cal R\cap\Lambda'(k)),\text{ if } x\in\Lambda'(k),\tag 5.33\\
10\text{ dist }(x,\Cal R\cap\Lambda''(k))&\geq\text{dist }(x,\partial\Lambda''(k))\geq\text{dist }
(x,\Cal R\cap\Lambda''(k)),\text{ if } x\in\Lambda''(k).\tag 5.34\endalign$$

Assume $x\in\Lambda'(k)$, ($x\in\Lambda''(k)$ works in the same way). 
We write $\Lambda'$ for $\Lambda'(k)$ for simplicity. Let $H_{\Lambda'\backslash\Cal R}$ be 
defined as in (3.14) and 
$$\Gamma=H-H_{\Lambda'\backslash\Cal R}.\tag 5.35$$
Then $$\phi(x)=\sum G_{\Lambda'\backslash\Cal R}(x,y)[\Gamma\phi](y),\tag 5.36$$
where $\phi$ is the eigenfunction of $H$ with eigenvalue $E$ ($\phi\neq 0$) and 
$$G_{\Lambda'\backslash\Cal R}(x,y)=(H_{\Lambda'\backslash\Cal R}-E)^{-1}(x,y)\tag 5.37$$

To estimate $G_{\Lambda'\backslash\Cal R}$, we use the resolvent equation. There are
two cases, $\Lambda\cap\Lambda'=\emptyset$ and $\Lambda\subset\Lambda'$. When
 $\Lambda\cap\Lambda'=\emptyset$, we cover $\Lambda'\backslash\Cal R$ with cubes 
 $\bar\Lambda$ and annulus $A=\bar\Lambda\backslash\{(\pm j, n)\}$, such that
 either (a) $\bar\Lambda\cap P=\emptyset$ or (b) $\bar\Lambda\cap P=(j',n')$
 and $$\text{dist }[(j',n'),\partial\bar\Lambda]\geq |j'|/4.\tag 5.38$$
 In case (a),
 $$\text{dist }(E, \sigma(H_{\bar\Lambda}))=\Cal O(1),\tag 5.39$$
 in case (b),
 $$\text{dist }(E, \sigma(H_{\bar\Lambda}))\geq \frac{\Cal O(\delta^2)}{[\min(|j|,|j'|)]^3},\tag 5.40$$
 similar to (5.26). In case (a), we have 
 $$|(H_{\bar\Lambda}-E)^{-1}(x,y)|\leq Ce^{-\kappa|x-y|},\,C,\kappa>0\tag 5.41$$
 for all $x$ and $y$. In case (b), we have 
 $$|(H_{\bar\Lambda}-E)^{-1}(x,y)|\leq Ce^{-\kappa|x-y|},\tag 5.42$$
 for all $x$ and $y$ such that $|x-y|\geq\sqrt {j'}$, $j'$ large by Neumann series.
 
 Iterating the resolvent equation using (5.41, 5.42, 5.27), we obtain 
 $$|G_{\Lambda'\backslash\Cal R}(x,y)|\leq e^{-\kappa'|x-y|},\,0<\kappa'<\kappa\tag 5.43$$
 for all $\Lambda'$, $x$, $y\in\Lambda'$, $|x-y|\geq k/10$. When $\Lambda'$ is such that 
 $\Lambda'\cap P=\{(j,n)\}$, then 
$$\align &\Vert (H_{\Lambda'\backslash\Cal R}-E)^{-1}\Vert=\Cal O(1)\tag 5.44\\
&|G_{\Lambda'\backslash\Cal R}(x,y)|\leq e^{-\kappa'|x-y|}\tag 5.45\endalign$$
for all $x$ and $y$.

When $\Lambda\subset\Lambda'$, we cover $\Lambda'\backslash\Cal R$ with $\Lambda$, annulus $A$ and
cubes $\bar\Lambda$ satisfying properties (a, b) as before. Using (5.41, 5.42, 5.27) in the 
resolvent equation, we obtain (5.43), assuming $|j|>K>K'$. (5.36, 5.43, 5.45) then
give (5.3). 

Using (5.28, 5.29) in place of (5.26, 5.27) in the resolvent iteration, we obtain (5.2) in the same way using (5.36) with $\Lambda$ replacing $\Cal R$.\hfill $\square$
\enddemo

\demo{Proof of Theorem 1}
We only need to verify that (H1, 3) are satisfied as Lemma 3.1 is available. This is where we need
$\delta$ to be small. Take $\Lambda=[-8, 8]^2$. From (2.3), $E\in\sigma(H_\Lambda)$ if and only if
$0\in\sigma(\tilde H_{\Lambda, E})$, where $\tilde H_{\Lambda, E}$ is defined similarly as in 
(2.3, 2.4) with $H_\Lambda$ replacing $H$ in (2.4). $\tilde H_{\Lambda, E}$ is a $5\times 5$ matrix.
Specifically we have 
$$
\align
\tilde H_{\Lambda, E} &=E-\delta\pmatrix 0&1&1&0&0\\ 1& 0&0&0&0\\
1& 0&0&0&0\\0& 0&0&0&0\\0& 0&0&0&0
\endpmatrix\\
& +\delta^2 \pmatrix 1&0&0&0&0\\ 0& 1/4&-1/2&0&0\\
0& -1/2&1/4&0&0\\0& 0&0&1/4&0\\0& 0&0&0&1/4
\endpmatrix+\Cal O(\delta^3)\\
& \ {\overset\text{def }\to=}E-A+\Cal O(\delta^3)\tag 5.46
\endalign
$$

It is easy to verify that the $5$ eigenvalues of $A$ satisfy 
$$\align &|\tilde\mu_1|,\,|\tilde\mu_2|\asymp\delta,\\
&|\tilde\mu_3|,\,|\tilde\mu_4|, |\tilde\mu_5|\asymp\delta^2,\, (0<\delta\ll 1)\tag 5.47\endalign$$
where $A\asymp B$ denotes $B/C\leq A\leq CB$ ($C>0$). So the eigenvalues $\mu$ of $H_\Lambda$
satisfy 
$$|\mu|\geq\Cal O(1)\delta^2.\tag 5.48$$
The eigenfunctions of $H_\Lambda$ satisfy 
$$\Vert \phi\Vert_{\ell^2(\partial\Lambda)}=\Cal O(\delta^4)\tag 5.49$$
by using the same arguments as in (5.6-5.11). (5.48, 5.49) imply that (H1, 3) are satisfied and 
in fact $0\notin\sigma(H_\Lambda)$ for all finite subsets $\Lambda\subset\Bbb Z^2$. So 
we reach the conclusion of Theorem 1.\hfill $\square$
\enddemo
\demo{Proof of Theorem 2} If $u(0)\in H^s(\Bbb T)$, ($s>0$), then $u(t)\in H^s(\Bbb T)$ for all $t$
and satisfies 
$$\Vert u(t)\Vert_{H^s}\leq C_s(1+|t|^s)\Vert u(0)\Vert_{H^s}.\tag 5.50$$
This holds generally for linear Schr\"odinger equation of the form (1.1) with real, smooth and bounded
$V$ (depending on $x$ and $t$), cf. Lemma 6.2 in \cite {B2}. 

Let $\phi=\phi(j,n)$ be the eigenfunctions of the Floquet Hamiltonian $H$ in (1.13) with eigenvalue $E$.
Then the Bloch waves $$\check\phi(x,t)=e^{iEt}\sum_{(j,n)\in\Bbb Z^2}\phi(j,n)e^{i(jx+nt)}\tag 5.51$$
provide a basis to expand the solution $u(t)$ with initial condition $u(0)$. More precisely for 
any given initial condition $u(0)\in L^2(\Bbb T)$, let $\hat u_0(j), \,j\in\Bbb Z$ be its 
Fourier coefficients. Identifying $\hat u_0\in\ell^2(\Bbb Z)$ with $\tilde u_0\in\ell^2(\Bbb Z^2)$ 
defined by 
$$\cases \tilde u_0(j,0)=\hat u_0(j)\\
\tilde u_0(j,n)=0,\,n\neq 0.\endcases\tag 5.52$$
We have 
$$u(x,t)=\sum_{\phi}(\tilde u_0,\phi)\check\phi(x,t).\tag 5.53$$
So $$\aligned \Vert u(t)\Vert^2_{H^s}&=\Vert \sum_{\phi}(\tilde u_0,\phi)\check\phi(t)\Vert^2_{H^s}\,(s>0)\\
&=\sum_j(1+|j|^{2s}) |\sum_{\phi}(\tilde u_0,\phi)(\check\phi(t), e_j)|^2\\
&=\sum_j(1+|j|^{2s}) |\sum_{\phi}\sum_k\hat u_0(k)\phi(k,0)(\check\phi(t), e_j)|^2,\endaligned\tag 5.54$$
where $$e_j=e^{-ijx},\tag 5.55$$and
$$\aligned (\check\phi(t), e_j)&=\int_0^{2\pi}\check\phi(t)e^{-ijx}dx\\
&=e^{iEt}\sum_{n\in\Bbb Z}\phi(j,n)e^{int}.\endaligned\tag 5.56$$

In view of the localization properties of $\phi$ in (1.14, 1.15), we decompose $\sum_k$ into
$\sum_{k,\,||k|-|j||\leq| j|/2}$ and $\sum_{k,\,||k|-|j||> |j|/2}$. We have 
$$\aligned \Vert u(t)\Vert^2_{H^s}&=\sum_j(1+|j|^{2s}) |\sum_{\phi}\sum_{k,\,||k|-|j||\leq| j/|2}\hat u_0(k)\phi(k,0)(\check\phi(t), e_j)|^2\\
&\quad +\sum_j(1+|j|^{2s}) |\sum_{\phi}\sum_{k,\,||k|-|j||>| j/|2}\hat u_0(k)\phi(k,0)(\check\phi(t), e_j)|^2\\
&\ {\overset\text{def }\to=}S_1+S_2.\endaligned\tag 5.57$$

Using (1.14, 1.15),
$$S_2\leq C\sum_j(1+|j|^{2s}) \sum_{k,\,||k|-|j||>| j/|2}e^{-\alpha||j|-|k||}\Vert u(0)\Vert_{H^s}\tag 5.58$$
for some $\alpha>0$ and we used $\Vert \hat u_0\Vert_{\infty}^2\leq\Vert \hat u_0\Vert_{2}^2\leq\Vert u(0)\Vert_{H^s}^2\,(s>0)$. So 
$$S_2\leq C_s^{(1)}\Vert u(0)\Vert_{H^s}^2.\tag 5.59$$

To estimate $S_1$, we notice that $||k|-|j||\leq| j/|2$. So 
$$\frac{|j|}{2}\leq|k|\leq\frac{3}{2}|j|\leq 2|j|.\tag 5.60$$
We make a dyadic expansion. Define 
$$u_0^{(\ell)}=\sum_{\frac{1}{4}\cdot 2^{\ell}\leq |k|\leq 4\cdot 2^{\ell}}\hat u_0(k)e^{ikx},\,\ell=0,\,1,...\tag 5.61$$
Then $$\aligned S_1&\leq C\sum_{\ell}(2^{\ell})^{2s}\sum_{{\frac{1}{2}}\cdot 2^{\ell}\leq |j|\leq 2\cdot 2^{\ell}}|\sum_{\phi}\sum_{k}\widehat {  u_0^{(\ell)}  }(k)\phi(k,0)(\check\phi(t), e_j)|^2\\
&\leq C\sum_\ell(2^\ell)^{2s}\Vert \sum_k\sum_\phi\widehat {  u_0^{(\ell)}  }(k)\phi(k,0)\check\phi(t)\Vert^2_2\\
&=C\sum_\ell(2^\ell)^{2s}\Vert \sum_\phi( u_0^{(\ell)},\phi)\check\phi(t)\Vert^2_2\\
&=C\sum_\ell(2^\ell)^{2s}\Vert u_0^{(\ell)}\Vert_2^2\\
&\leq C_s^{(2)}\Vert u(0)\Vert_{H^s}^2,\endaligned\tag 5.62$$
where we used (5.60, 5.61), $\ell^2$ norm conservation. The last line follows by standard considerations
using dyadic expansion. Using (5.59, 5.62) in (5.57), we obtain (1.16).\hfill $\square$
\enddemo

\enddocument
\end